\documentclass[]{amsart}
\usepackage{amsmath, amsthm, amscd, amsfonts, amssymb, graphicx, color}
\usepackage[bookmarksnumbered, colorlinks, plainpages]{hyperref}

\theoremstyle{definition}

\theoremstyle{remark}

\numberwithin{equation}{section}

\begin{document}
\setcounter{page}{1}
\begin{center}
{\bf ARENS REGULARITY AND MODULE ARENS REGULARITY OF MODULE ACTIONS}
\end{center}
\address{}
\address{}
\title[]{}
\dedicatory{}
\title[]{}
\author[]{KAZEM HAGHNEJAD AZAR    }

\address{}

\address{}

\dedicatory{}



\begin{abstract}
 In this paper,  we  extend some problems  from Arens regularity and module Arens regularity of Banach algebras to module actions.
 \end{abstract} \maketitle

\begin{center}
{\bf  1.Introduction and Preliminaries}
\end{center}
\vspace{0.5cm}
\noindent In [13], for Banach algebras $A$ and $\mathcal{O}$, authors extended the concept of Arens regularity of Banach algebra $A$ to the case that there is an $\mathcal{O}$-module structure of $A$ which is called module Arens regularity of $A$ as $\mathcal{O}$-module. In this note, we study this problem for a right module action $\pi_r :B\times A\rightarrow B$ where $B$ is a Banach $A$-bimodule and we extend some problems from Arens regularity of Banach algebras to the left and right module actions.\\
Let $X,Y,Z$ be normed spaces and $m:X\times Y\rightarrow Z$ be a bounded bilinear mapping. Arens  offers two natural extensions $m^{***}$ and $m^{t***t}$ of $m$ from $X^{**}\times Y^{**}$ into $Z^{**}$ as follows:\\
 \noindent1. $m^*:Z^*\times X\rightarrow Y^*$,~~~~~given by~~~$\langle  m^*(z^\prime,x),y\rangle  =\langle  z^\prime, m(x,y)\rangle  $ ~where $x\in X$, $y\in Y$, $z^\prime\in Z^*$,\\
 2. $m^{**}:Y^{**}\times Z^{*}\rightarrow X^*$,~~given by $\langle  m^{**}(y^{\prime\prime},z^\prime),x\rangle  =\langle  y^{\prime\prime},m^*(z^\prime,x)\rangle  $ ~where $x\in X$, $y^{\prime\prime}\in Y^{**}$, $z^\prime\in Z^*$,\\
3. $m^{***}:X^{**}\times Y^{**}\rightarrow Z^{**}$,~ given by~ ~ ~$\langle  m^{***}(x^{\prime\prime},y^{\prime\prime}),z^\prime\rangle  $ $=\langle  x^{\prime\prime},m^{**}(y^{\prime\prime},z^\prime)\rangle  $ \\~where ~$x^{\prime\prime}\in X^{**}$, $y^{\prime\prime}\in Y^{**}$, $z^\prime\in Z^*$.\\
The mapping $m^{***}$ is the unique extension of $m$ such that $x^{\prime\prime}\rightarrow m^{***}(x^{\prime\prime},y^{\prime\prime})$ from $X^{**}$ into $Z^{**}$ is $weak^*-to-weak^*$ continuous for every $y^{\prime\prime}\in Y^{**}$, but the mapping $y^{\prime\prime}\rightarrow m^{***}(x^{\prime\prime},y^{\prime\prime})$ is not in general $weak^*-to-weak^*$ continuous from $Y^{**}$ into $Z^{**}$ unless $x^{\prime\prime}\in X$. Hence the first topological center of $m$ may  be defined as following
$$Z_1(m)=\{x^{\prime\prime}\in X^{**}:~~y^{\prime\prime}\rightarrow m^{***}(x^{\prime\prime},y^{\prime\prime})~~is~~weak^*-to-weak^*-continuous\}.$$
Let  $m^t:Y\times X\rightarrow Z$ be the transpose of $m$ defined by $m^t(y,x)=m(x,y)$ for every $x\in X$ and $y\in Y$. Then $m^t$ is a continuous bilinear map from $Y\times X$ to $Z$, and so it may be extended as above to $m^{t***}:Y^{**}\times X^{**}\rightarrow Z^{**}$.
 The mapping $m^{t***t}:X^{**}\times Y^{**}\rightarrow Z^{**}$ in general is not equal to $m^{***}$, see [1], if $m^{***}=m^{t***t}$, then $m$ is called Arens regular. The mapping $y^{\prime\prime}\rightarrow m^{t***t}(x^{\prime\prime},y^{\prime\prime})$ is $weak^*-to-weak^*$ continuous for every $y^{\prime\prime}\in Y^{**}$, but the mapping $x^{\prime\prime}\rightarrow m^{t***t}(x^{\prime\prime},y^{\prime\prime})$ from $X^{**}$ into $Z^{**}$ is not in general  $weak^*-to-weak^*$ continuous for every $y^{\prime\prime}\in Y^{**}$. So we define the second topological center of $m$ as
$$Z_2(m)=\{y^{\prime\prime}\in Y^{**}:~~x^{\prime\prime}\rightarrow m^{t***t}(x^{\prime\prime},y^{\prime\prime})~~is~~weak^*-to-weak^*-continuous\}.$$
It is clear that $m$ is Arens regular if and only if $Z_1(m)=X^{**}$ or $Z_2(m)=Y^{**}$. Arens regularity of $m$ is equivalent to the following
$$\lim_i\lim_j\langle  z^\prime,m(x_i,y_j)\rangle  =\lim_j\lim_i\langle  z^\prime,m(x_i,y_j)\rangle  ,$$
whenever both limits exist for all bounded sequences $(x_i)_i\subseteq X$ , $(y_i)_i\subseteq Y$ and $z^\prime\in Z^*$, see [5, 12].\\
 The regularity of a Banach algebra $A$ is defined to be the regularity of its algebra multiplication when considered as a bilinear mapping such as $m$. Let $a^{\prime\prime}$ and $b^{\prime\prime}$ be elements of $A^{**}$, the second dual of $A$. By $Goldstin^,s$ Theorem [6, P.424-425], there are nets $(a_{\alpha})_{\alpha}$ and $(b_{\beta})_{\beta}$ in $A$ such that $a^{\prime\prime}=weak^*-\lim_{\alpha}a_{\alpha}$ ~and~  $b^{\prime\prime}=weak^*-\lim_{\beta}b_{\beta}$. So it is easy to see that for all $a^\prime\in A^*$, we have
$$\lim_{\alpha}\lim_{\beta}\langle  a^\prime,m(a_{\alpha},b_{\beta})\rangle  =\langle  a^{\prime\prime}b^{\prime\prime},a^\prime\rangle  $$ and
$$\lim_{\beta}\lim_{\alpha}\langle  a^\prime,m(a_{\alpha},b_{\beta})\rangle  =\langle  a^{\prime\prime}ob^{\prime\prime},a^\prime\rangle  ,$$
where $a^{\prime\prime}b^{\prime\prime}$ and $a^{\prime\prime}ob^{\prime\prime}$ are the first and second Arens products of $A^{**}$, respectively, see [10, 12].\\
The mapping $m$ is left strongly Arens irregular if $Z_1(m)=X$ and $m$ is right strongly Arens irregular if $Z_2(m)=Y$.\\

\begin{center}
\textbf{{ 2. The topological centers of module actions}}\\
\end{center}

\vspace{0.5cm}

 Let $B$ be a Banach $A-bimodule$, and let\\
$$\pi_\ell:~A\times B\rightarrow B~~~and~~~\pi_r:~B\times A\rightarrow B.$$
be the left and right module actions of $A$ on $B$. Then $B^{**}$ is a Banach $A^{**}-bimodule$ with module actions
$$\pi_\ell^{***}:~A^{**}\times B^{**}\rightarrow B^{**}~~~and~~~\pi_r^{***}:~B^{**}\times A^{**}\rightarrow B^{**}.$$
Similarly, $B^{**}$ is a Banach $A^{**}-bimodule$ with module actions\\
$$\pi_\ell^{t***t}:~A^{**}\times B^{**}\rightarrow B^{**}~~~and~~~\pi_r^{t***t}:~B^{**}\times A^{**}\rightarrow B^{**}.$$
We may therefore define the topological centers of the right and left module actions of $A$ on $B$ as follows:
$$Z_{A^{**}}(B^{**})=Z(\pi_r)=\{b^{\prime\prime}\in B^{**}:~the~map~~a^{\prime\prime}\rightarrow \pi_r^{***}(b^{\prime\prime}, a^{\prime\prime})~:~A^{**}\rightarrow B^{**}$$$$~is~~~weak^*-to-weak^*~continuous\}$$
$$Z_{B^{**}}(A^{**})=Z(\pi_\ell)=\{a^{\prime\prime}\in A^{**}:~the~map~~b^{\prime\prime}\rightarrow \pi_\ell^{***}(a^{\prime\prime}, b^{\prime\prime})~:~B^{**}\rightarrow B^{**}$$$$~is~~~weak^*-to-weak^*~continuous\}$$
$$Z_{A^{**}}^t(B^{**})=Z(\pi_\ell^t)=\{b^{\prime\prime}\in B^{**}:~the~map~~a^{\prime\prime}\rightarrow \pi_\ell^{t***}(b^{\prime\prime}, a^{\prime\prime})~:~A^{**}\rightarrow B^{**}$$$$~is~~~weak^*-to-weak^*~continuous\}$$
$$Z_{B^{**}}^t(A^{**})=Z(\pi_r^t)=\{a^{\prime\prime}\in A^{**}:~the~map~~b^{\prime\prime}\rightarrow \pi_r^{t***}(a^{\prime\prime}, b^{\prime\prime})~:~B^{**}\rightarrow B^{**}$$$$~is~~~weak^*-to-weak^*~continuous\}$$
We note also that if $B$ is a left(resp. right) Banach $A-module$ and $\pi_\ell:~A\times B\rightarrow B$~(resp. $\pi_r:~B\times A\rightarrow B$) is left (resp. right) module action of $A$ on $B$, then $B^*$ is a right (resp. left) Banach $A-module$. \\
We write $ab=\pi_\ell(a,b)$, $ba=\pi_r(b,a)$, $\pi_\ell(a_1a_2,b)=\pi_\ell(a_1,a_2b)$,  $\pi_r(b,a_1a_2)=\pi_r(ba_1,a_2)$,~\\
$\pi_\ell^*(a_1b^\prime, a_2)=\pi_\ell^*(b^\prime, a_2a_1)$,~
$\pi_r^*(b^\prime a, b)=\pi_r^*(b^\prime, ab)$,~ for all $a_1,a_2, a\in A$, $b\in B$ and  $b^\prime\in B^*$
when there is no confusion.\\
A functional $a^\prime$ in $A^*$ is said to be $wap$ (weakly almost
 periodic) on $A$ if the mapping $a\rightarrow a^\prime a$ from $A$ into
 $A^{*}$ is weakly compact.  In  [12], Pym showed that  this definition to the equivalent following condition\\
 For any two net $(a_{\alpha})_{\alpha}$ and $(b_{\beta})_{\beta}$
 in $\{a\in A:~\parallel a\parallel\leq 1\}$, we have\\
$$\\lim_{\alpha}\\lim_{\beta}\langle a^\prime,a_{\alpha}b_{\beta}\rangle=\\lim_{\beta}\\lim_{\alpha}\langle a^\prime,a_{\alpha}b_{\beta}\rangle,$$
whenever both iterated limits exist. The collection of all $wap$
functionals on $A$ is denoted by $wap(A)$. Also we have
$a^{\prime}\in wap(A)$ if and only if $\langle a^{\prime\prime}b^{\prime\prime},a^\prime\rangle=\langle a^{\prime\prime}ob^{\prime\prime},a^\prime\rangle$ for every $a^{\prime\prime},~b^{\prime\prime} \in
A^{**}$. \\
Let $B$ be a Banach left $A-module$. Then, $b^\prime\in B^*$ is said to be left weakly almost periodic functional if the set $\{\pi^*_\ell(b^\prime,a):~a\in A,~\parallel a\parallel\leq 1\}$ is relatively weakly compact. We denote by $wap_\ell(B)$ the closed subspace of $B^*$ consisting of all the left weakly almost periodic functionals in $B^*$.\\
The definition of the right weakly almost periodic functional ($=wap_r(B)$) is the same.\\
By  [12], $b^\prime\in wap_\ell(B)$ is equivalent to the following $$\langle \pi_\ell^{***}(a^{\prime\prime},b^{\prime\prime}),b^\prime\rangle=
\langle \pi_\ell^{t***t}(a^{\prime\prime},b^{\prime\prime}),b^\prime\rangle$$
for all $a^{\prime\prime}\in A^{**}$ and $b^{\prime\prime}\in B^{**}$.
Thus, we can write \\
$$wap_\ell(B)=\{ b^\prime\in B^*:~\langle \pi_\ell^{***}(a^{\prime\prime},b^{\prime\prime}),b^\prime\rangle=
\langle \pi_\ell^{t***t}(a^{\prime\prime},b^{\prime\prime}),b^\prime\rangle~~$$$$for~~all~~a^{\prime\prime}\in A^{**},~b^{\prime\prime}\in B^{**}\}.$$\\\\

\noindent {\it{\bf Theorem 2-1.}} Suppose that $B$ is a left Banach $A-module$. Then  the following assertions are equivalents.
\begin{enumerate}
\item The mapping $a\rightarrow \pi^*_\ell (b^\prime , a)$ from $A$ into $B^*$ is weakly compact.
\item ${Z}^\ell_{B^{**}}(A^{**})=A^{**}$.

\item There are a subset $E$ of $B^*$ with $\overline{lin}E=B^*$ such that for each sequence $(a_n)_n\subseteq A$ and $(b_m)_m\subseteq B$ and each $b^\prime\in B^*$, we have
    $$\lim_m\lim_n \langle  b^\prime,a_nb_m\rangle=\lim_n\lim_m \langle  b^\prime,a_nb_m\rangle,$$
    whenever both the iterated limits exist.\\
\item Suppose that $a^{\prime \prime}\in A^{**}$ and $(a_\alpha)_\alpha \subseteq A$ such that $a_\alpha \stackrel{w^*} {\rightarrow}a^{\prime \prime}$. Then we have
     $$\pi^*_\ell (b^\prime , a_\alpha)\stackrel{w} {\rightarrow}\pi^{t**t}_\ell (b^\prime , a^{\prime \prime}),$$ for each $b^\prime\in B^*$.
\end{enumerate}

\begin{proof}
$(1)\Rightarrow (2)$\\
Suppose that $b^\prime\in B^*$. Take $T(a)=\pi^*_\ell (b^\prime , a)$ where $a\in A$. By easy calculation, we have $T^{**}(a^{\prime \prime})=\pi_\ell^{****} (b^\prime , a^{\prime \prime})$ for each $a^{\prime \prime}\in A^{**}$. Now let $T$ be a weakly compact mapping. Then by using Theorem VI 4.2 and VI 4.8,  from [6], we have $\pi_\ell^{****} (b^\prime , a^{\prime \prime})\in B^*$ for each  $a^{\prime \prime}\in A^{**}$. Suppose that $(b_{\alpha}^{\prime\prime})_{\alpha}\subseteq B^{**}$ such that  $b^{\prime\prime}_{\alpha} \stackrel{w^*} {\rightarrow}b^{\prime\prime}$ on $B^{**}$. Then for every $a^{\prime \prime}\in A^{**}$, we have
$$ \langle  \pi_\ell^{***} (a^{\prime \prime},b^{\prime\prime}_{\alpha}), b^\prime\rangle=
 \langle  a^{\prime \prime},\pi_\ell^{**} (b^{\prime\prime}_{\alpha}, b^\prime )\rangle=
 \langle  \pi_\ell^{*****} (b^{\prime\prime}_{\alpha}, b^\prime ),a^{\prime \prime}\rangle$$$$=
 \langle  b^{\prime\prime}_{\alpha},\pi_\ell^{****} ( b^\prime ,a^{\prime \prime})\rangle
\rightarrow  \langle  b^{\prime\prime},\pi_\ell^{****} (b^\prime ,a^{\prime \prime})\rangle= \langle  \pi_\ell^{***} (a^{\prime \prime},b^{\prime\prime}), b^\prime\rangle.$$
It follows that $\pi_\ell^{***} (a^{\prime \prime},b^{\prime\prime}_{\alpha})\stackrel{w^*} {\rightarrow}\pi_\ell^{***} (a^{\prime \prime},b^{\prime\prime})$, and so $a^{\prime \prime}\in {Z}^\ell_{B^{**}}(A^{**})$.\\
$(2)\Rightarrow (1)$\\
Let ${Z}^\ell_{B^{**}}(A^{**})=A^{**}$. Suppose that $(b_{\alpha}^{\prime\prime})_{\alpha}\subseteq B^{**}$ such that  $b^{\prime\prime}_{\alpha} \stackrel{w^*} {\rightarrow}b^{\prime\prime}$ in $B^{**}$. Then for every $a^{\prime \prime}\in A^{**}$, we have $\pi_\ell^{***} (a^{\prime \prime},b^{\prime\prime}_{\alpha})\stackrel{w^*} {\rightarrow}\pi_\ell^{***} (a^{\prime \prime},b^{\prime\prime})$. It follows that
$$ \langle  T^{**}(a^{\prime \prime}),b_{\alpha}^{\prime\prime}\rangle\rightarrow  \langle  T^{**}(a^{\prime \prime}),b^{\prime\prime}\rangle,$$
for each $a^{\prime \prime}\in A^{**}$. Consequently, $T^{**}(a^{\prime \prime})\in B^*$ for each $a^{\prime \prime}\in A^{**}$, and so $T^{**}(A^{\prime \prime})\subseteq B^*$. By another using Theorem VI 4.2 and VI 4.8,  from [6], we conclude that the mapping $a\rightarrow\pi_\ell (b^\prime , a)$ from $A$ into $B^*$ is weakly compact.\\
$(2)\Rightarrow (3)$\\
By  definition of ${Z}^\ell_{B^{**}}(A^{**})$, since ${Z}^\ell_{B^{**}}(A^{**})=A^{**}$, proof  hold.\\
$(3)\Rightarrow (1)$\\
Proof is similar to Theorem 2.6.17 from [5].\\
$(1)\Rightarrow (4)$\\
Let $a^{\prime \prime}\in A^{**}$ and $(a_\alpha)_\alpha \subseteq A$ such that $a_\alpha \stackrel{w^*} {\rightarrow}a^{\prime \prime}$. Then for each $b^{\prime \prime}\in B^{**}$, we have
$$\lim_\alpha  \langle  b^{\prime \prime}, \pi^*_\ell (b^\prime , a_\alpha)\rangle=\lim_\alpha  \langle  \pi^{**}_\ell (b^{\prime \prime}, b^\prime ), a_\alpha)\rangle= \langle  \pi_\ell^{***} (a^{\prime \prime},b^{\prime\prime}), b^\prime\rangle$$$$=
 \langle  \pi_\ell^{t***t} (a^{\prime \prime},b^{\prime\prime}), b^\prime\rangle= \langle  b^{\prime \prime},\pi_\ell^{t**} (a^{\prime \prime},b^{\prime})\rangle.$$
It follows that  $\pi^*_\ell (b^\prime , a_\alpha)\stackrel{w} {\rightarrow}\pi^{t**t}_\ell (b^\prime , a^{\prime \prime})$, so this completes the proof.\\
$(4)\Rightarrow (2)$\\
Let $b^\prime\in B^*$ and suppose that $a^{\prime \prime}\in A^{**}$ and $b^{\prime \prime}\in B^{**}$. Let $(a_\alpha)_\alpha \subseteq A$ such that $a_\alpha \stackrel{w^*} {\rightarrow}a^{\prime \prime}$. Since
     $$\pi^*_\ell (b^\prime , a_\alpha)\stackrel{w} {\rightarrow}\pi^{t**t}_\ell (b^\prime , a^{\prime \prime}),$$ for each $b^\prime\in B^*$, we have the following equality
$$ \langle  \pi_\ell^{***} (a^{\prime \prime},b^{\prime\prime}), b^\prime\rangle=  \langle  a^{\prime \prime},\pi_\ell^{**} (b^{\prime\prime}, b^\prime)\rangle=\lim_\alpha \langle  \pi_\ell^{**} (b^{\prime\prime}, b^\prime),a_\alpha\rangle$$$$=
 \lim_\alpha \langle  b^{\prime\prime}, \pi_\ell^{*} (b^\prime,a_\alpha )\rangle=  \langle  b^{\prime\prime},\pi^{t**t}_\ell (b^\prime , a^{\prime \prime})$$$$= \langle  \pi_\ell^{t***t} (a^{\prime \prime},b^{\prime\prime}), b^\prime\rangle.$$
It follows that $b^\prime \in wap_\ell(B)$, and so  ${Z}^\ell_{B^{**}}(A^{**})=A^{**}$.
\end{proof}
\vspace{0.5cm}

\noindent {\it{\bf Corollary 2-2.}} Suppose that $B$ is a left Banach $A-module$. Then $B^*A^{**}\subseteq B^*$ if and only if ${Z}^\ell_{B^{**}}(A^{**})=A^{**}$.\\

\vspace{0.5cm}

\noindent {\it{\bf Example 2-3.}} Suppose that  $G$ is a locally compact group. In the preceding corollary, take $A=B=c_0(G)$. Therefore we conclude that
${Z}^\ell_1(\ell^1(G)^{**})=\ell^1(G)^{**},$ see [5, Example 2.6.22(iii)].\\

\vspace{0.5cm}

\noindent {\it{\bf Theorem 2-4.}} Suppose that $B$ is a right Banach $A-module$. Then  the following assertions are equivalents.
\begin{enumerate}

\item ${Z}^\ell_{A^{**}}(B^{**})=B^{**}$.

\item The mapping $b\rightarrow \pi_r^{*}(b^\prime ,b)$ from $B$ into $A^*$ is weakly compact.

\item There are a subset $E$ of $B^*$ with $\overline{lin}E=B^*$ such that for each sequence $(a_n)_n\subseteq A$ and $(b_m)_m\subseteq B$ and each $b^\prime\in B^*$,
    $$\lim_m\lim_n \langle  b^\prime,b_ma_n\rangle=\lim_n\lim_m \langle  b^\prime,b_ma_n\rangle,$$
    whenever both the iterated limits exist.\\
\item Suppose that $b^{\prime \prime}\in B^{**}$ and $(b_\alpha)_\alpha\subseteq B$ such that $b_{\alpha} \stackrel{w^*} {\rightarrow}b^{\prime\prime}$ on $B^{**}$. Then
 $$\pi_r^*(b^\prime,b_\alpha)\stackrel{w} {\rightarrow} \pi_r^{t**t}(b^\prime,b^{\prime\prime}),$$
 for all $b^\prime \in B^*$.

\end{enumerate}
\begin{proof}
Proof is similar to Theorem 2-1.\end{proof}

\vspace{0.5cm}

\noindent {\it{\bf Corollary 2-5.}} Suppose that $B$ is a right Banach $A-module$. Then $B^*B^{**}\subseteq A^*$ if and only if ${Z}^\ell_{A^{**}}(B^{**})=B^{**}$.\\

\vspace{0.5cm}

 \noindent {\it{\bf Corollary 2-6.}} In the preceding corollary, if we take $B=A$, we obtain Lemma 3.1 (i) from [10] and in the Theorems 2-1, if we take $B=A$,  we obtain Theorem 2.6.17 from [5].\\

\vspace{0.5cm}

\noindent {\it{\bf Definition 2-7.}} The bilinear mapping $\pi_r:B\times A\rightarrow B$ is called module Arens regular if satisfies in the following conditions:\\
i) the mapping $T_{b^\prime}:b\rightarrow b^\prime b$ from $B$ into $A^*$ is weakly compact for any $b^\prime\in B^*$ for which the mapping $T_{b^\prime}$ is $A-$module homomorphism.\\
ii) $T_{b^\prime a}$ from $B$ into $A^*$ is $A-$module homomorphism when $T_{b^\prime}$ is $A-$module homomorphism $b^\prime\in B^*$ and $a\in A$.\\

\vspace{0.5cm}

Suppose that $B$ is a  Banach $A-bimodule$. We assume that $J$ is the closed right ideal of $A$ generated by elements of the form $a_1ba_2-ba_1a_2$ for all $a_1,a_2\in A$, $b\in B$.\\

\vspace{0.5cm}

 \noindent {\it{\bf Theorem 2-8.}} Take $J$ is defined as above. Then $\pi_r$ is a module Arens regular if and only if
 $\pi_r^{***}(b^{\prime \prime},a^{\prime \prime})-\pi_r^{t***t}(b^{\prime \prime},a^{\prime \prime})\in J^{\perp\perp}$ for $a^{\prime \prime}\in A^{**}$ and $b^{\prime \prime}\in B^{**}$.
\begin{proof} The mapping $T_{b^\prime}:b\rightarrow b^\prime b$ from $B$ into $A^*$ is an $A$-module homomorphism if and only if $T_{b^\prime}(ab)=aT_{b^\prime}(b)$ for all $b\in B$ and $a\in A$,
$$\Leftrightarrow b^\prime ab=ab^\prime b,$$
$$\Leftrightarrow  <b^\prime, abx-bxa>=0~~~for~~all~~x\in A.$$
In the above statements, if we replace $T_{b^\prime y}$ with $T_{b^\prime}$ where $y\in A$, then the last equality is equivalent with $b^\prime\in J^\perp$. Consequently, the bilinear mapping $\pi_r$ is module Arens regular if and only $T_{b^\prime}$  is weakly compact for any $b^\prime\in J^\perp$. By using Theorem 2.4, $T_{b^\prime}$  is weakly compact if and only if for every $a^{\prime \prime}\in A^{**}$ and $b^{\prime \prime}\in B^{**}$ we have $<\pi_r^{***}(b^{\prime \prime},a^{\prime \prime}),b^\prime>=<\pi_r^{t***t}(b^{\prime \prime},a^{\prime \prime}),b^\prime>$, and so $\pi_r^{***}(b^{\prime \prime},a^{\prime \prime})-\pi_r^{t***t}(b^{\prime \prime},a^{\prime \prime})\in J^{\perp\perp}$. This complete the proof.\end{proof}

\vspace{0.5cm}

\noindent {\it{\bf Theorem 2-9.}} Suppose that $B$ is a right Banach $A-module$. Then  the following assertions are equivalents.
\begin{enumerate}

\item $\pi_r$ is module Arens regular.

\item the mapping $T_{b^\prime}:b\rightarrow b^\prime b$ from $B$ into $A^*$ is weakly compact and  $T_{b^\prime a}$ from $B$ into $A^*$ is $A-$module homomorphism  for which the mapping $T_{b^\prime}$ is $A-$module homomorphism
 whenever $b^\prime\in B^*$ and $a\in A$.\\
\item There are a subset $E$ of $B^*$ with $\overline{lin}E=B^*$ such that for each sequence $(a_n)_n\subseteq A$ and $(b_m)_m\subseteq B$ and each $b^\prime\in J^\perp$,
    $$\lim_m\lim_n \langle  b^\prime,b_ma_n\rangle=\lim_n\lim_m \langle  b^\prime,b_ma_n\rangle,$$
    whenever both the iterated limits exist.\\
\item the mapping $b^{\prime \prime}\rightarrow \pi_r^{***}(b^{\prime \prime},a^{\prime \prime})$ is $\sigma(B^{**},J^\perp)-$continuous for every $a^{\prime \prime}\in A^{**}$.

\end{enumerate}
\begin{proof}
Proof is similar to the preceding theorem and Theorem 2.4.\end{proof}

\vspace{0.5cm}

\bibliographystyle{amsplain}

\it{Department of Mathematics, Amirkabir University of Technology, Tehran, Iran\\
{\it Email address:} haghnejad@aut.ac.ir\\\\
Department of Mathematics, Amirkabir University of Technology, Tehran, Iran\\
{\it Email address:} riazi@aut.ac.ir}

\end{document}